\documentstyle[11pt,amssymb,amsfonts]{article}

\begin{document}
\newcommand{\p}{\parallel }
\makeatletter \makeatother
\newtheorem{th}{Theorem}[section]
\newtheorem{lem}{Lemma}[section]
\newtheorem{de}{Definition}[section]
\newtheorem{rem}{Remark}[section]
\newtheorem{cor}{Corollary}[section]
\renewcommand{\theequation}{\thesection.\arabic {equation}}

\title{{\bf The Greiner's approach to heat kernel asymptotics and the variation formulas for the equivariant
Ray-Singer metric }
}
\author{ Yong Wang }

\date{}
\maketitle

\begin{abstract} In this paper,
using the Greiner's approach to heat kernel asymptotics, we give new
proofs of the equivariant Gauss-Bonnet-Chern formula and the
variation formulas for the equivariant Ray-Singer metric, which are
originally due to J. M. Bismut and W. Zhang.
\\

\\
\end{abstract}

\section{Introduction}
    \quad The first success of proving the Atiyah-Singer index theorem directly by heat kernel method
was achieved by Patodi \cite{Pa}, who carried out the "fantastic
cancellation" (cf. \cite{MS}) for the Laplace operators and for the
first time proved a local version of the Gauss-Bonnet-Chern theorem.
Later on
 several different direct heat kernel proofs of the Atiyah-Singer index theorem for
Dirac operators have appeared independently: Bismut \cite{Bi}, Getzler
\cite{Ge1}, \cite{Ge2},  Yu \cite{Yu} and  Ponge \cite{Po}. All the proofs have their own
advantages.\\
\indent The Atiyah-Bott-Segal-Singer index formula is a
generalization with group action of the Atiyah-Singer index theorem.
In \cite{BV}, Berline and Vergne gave a heat kernel proof of the
Atiyah-Bott-Segal-Singer index formula. In \cite{LYZ}, Lafferty, Yu and
Zhang presented a very simple and direct geometric proof for
equivariant index of the Dirac operator. In \cite{PW}, Ponge and Wang
gave a different proof of the equivariant index formula by the
Greiner's approach to the heat kernel asymptotics. Based on the method in \cite{LYZ}, Zhou gave a direct geometric proof of the equivariant Gauss-Bonnet-Chern formula in \cite{Zh}. THE FIRST PURPOSE of this paper is to give another proof of the equivariant Gauss-Bonnet-Chern formula
 by the Greiner's approach to the heat kernel asymptotics.  \\
   \indent In \cite{BZ1}, Bismut and Zhang extended the famous Cheeger-M\"{u}ller theorem to the case where the metric on the auxiliary bundle is not flat
   and they proved anomaly formulas for Ray-Singer metrics. In \cite{BZ2}, Bismut and Zhang further extended their generalized Cheeger-M\"{u}ller theorem to
   the equivariant case and gave anomaly formulas for equivariant Ray-Singer metrics. In \cite{We}, Weiss gave a new and detailed proof of the variation formulas for the
equivariant Ray-Singer metric due to Bismut-Zhang by the method in \cite{BV}. The proof of Weiss and Berline-Vergne lifted the operators to the principle
bundle. THE SECOND PURPOSE of this paper is to give another proof of anomaly formulas for the
equivariant Ray-Singer metric due to Bismut-Zhang
 by the Greiner's approach to the heat kernel asymptotics. In our proof, we do not need to
lift operators to the principle bundle and use a trick due to Chern and Hu in \cite{CH}.
Using the approach to Greiner's heat kernel asymptotics to give a new proof of the
variation formulas for the equivariant Ray-Singer metric has two advantages. One is to use
the Volterra pseudodifferential calculus to get the heat kernel asymtotic expansion
instead of heat equation discussions as in \cite{BGV}, \cite{BZ2}. In \cite{LYZ}, for
proving the equivariant local index theorem, transformed formulas between normal
and tubular coordinates are needed. The other advantage is that the transformed
formulas between normal and tubular coordinates are the consequence of a standard
change of variable formula for pseudodifferential symbols.\\
    \indent
    This paper is organized as follows: In Section 2, we give another proof of the equivariant Gauss-Bonnet-Chern formula.
  In Section 3, we give another proof of the variation formulas for the
equivariant Ray-Singer metric.\\

\section{The equivariant Gauss-Bonnet-Chern formula}

 \quad Let $M$ be a closed even dimensional $n$ oriented Riemannian manifold and $\phi$ be an isometry on $M$ preserving the orientation.
 Then $\phi$ induces a map
 $$\widetilde{\phi}={\phi^{-1}}^*:\wedge T^*_xM\rightarrow \wedge T^*_{\phi x}M$$
  on the exterior algebra bundle $\wedge T^*M$. Let $d$ denote the exterior differential operator and $\delta$
 be its adjoint operator and $D=d+\delta$ be the de-Rham Hodge operator. Let $D^+=D|_{\wedge^{\rm even} T^*M}$ and $D^-=D|_{\wedge^{\rm odd} T^*M}$. Then $\widetilde{\phi}D=D\widetilde{\phi}$ and we define the equivariant
 index $${\rm Ind}_\phi(D)={\rm Tr}(\widetilde{\phi}|_{{\rm ker}D^+})-{\rm Tr}(\widetilde{\phi}|_{{\rm ker}D^-}).\eqno(2.1)$$
\indent We recall the Greiner's approach to heat kernel
asymptotics as in
\cite{BeGS}, p.362, \cite{Po}, p.216. Define
 the operator given by
 $$(Q_0u)(x,s)=\int_{0}^{\infty}e^{-sD^2}[u(x,t-s)]dt,~~u\in \Gamma_c(M\times {\mathbb{R}},\wedge T^*M),\eqno(2.2)$$
 maps continuously $u$ to $D'(M\times {\mathbb{R}},\wedge T^*M))$ which is the dual space of $\Gamma_c(M\times {\mathbb{R}},\wedge T^*M)).$ We have
 $$(D^2+\frac{\partial}{\partial t})Q_0u=Q_0(D^2+\frac{\partial}{\partial t})u=u,~~~u\in \Gamma_c(M\times {\mathbb{R}},\wedge T^*M)). \eqno(2.3)$$
Let $(D^2+\frac{\partial}{\partial t})^{-1}$ be the Volterra inverse of $D^2+\frac{\partial}{\partial t}$ as in \cite{BeGS}, p.362. Then
$$(D^2+\frac{\partial}{\partial t})Q=I-R_1;~~Q(D^2+\frac{\partial}{\partial t})=1-R_2,\eqno(2.4)$$
where $R_1,R_2$ are smooth operators. Let
$$(Q_0u)(x,t)=\int_{M\times {\mathbb{R}}}K_{Q_0}(x,y,t-s)u(y,s)dyds,\eqno(2.5)$$
and $k_t(x,y)$ is the heat kernel of $e^{-tD^2}$. We get
$$K_{Q_0}(x,y,t)=k_t(x,y)~ {\rm when}~ t>0,~~ {\rm when }~ t<0,~ K_{Q_0}(x,y,t)=0.\eqno(2.6)$$

\noindent{\bf Definition 2.1}~~{\it The operator $P$ is called the Volterra $\Psi DO$ if
 (i) $P$ has the Volterra property,i.e. it has a distribution
kernel of the form $K_P(x,y, t-s)$ where $K_P(x, y, t)$ vanishes on the region $t < 0.$}\\
(ii) {\it The parabolic homogeneity of the heat operator $P +
\frac{\partial}{\partial t}$, i.e. the homogeneity with respect to
the dilations of ${\mathbb{R}}^n\times{\mathbb{R}} ^1$ is given by}
 $$\lambda\cdot (\xi,\tau)=(\lambda\xi,\lambda^2\tau),~~~~~~~(\xi,\tau)\in {\mathbb{R}}^n\times {\mathbb{R}}^1,~~\lambda\neq 0.\eqno(2.7)$$\\

\indent In the sequel for $g\in  \textsl{S} ({\mathbb{R}}^{n+1})$
and $\lambda\neq 0$, we let $g_{\lambda}$ be the tempered
distribution defined by
$$\left<g_\lambda(\xi,\tau),u(\xi,\tau)\right>=|\lambda|^{-(n+2)}\left<g_\lambda(\xi,\tau),u(\lambda^{-1}\xi,\lambda^{-2}\tau)\right>,~~u\in
{\textsl{S}} ({\mathbb{R}}^{n+1}).\eqno(2.8)$$\\

\noindent{\bf Definition 2.2}~~A distribution $ g\in  \textsl{S}
({\mathbb{R}}^{n+1})$ is parabolic homogeneous of degree $m,~ m \in{\mathbb{Z}},
$ if for
any $\lambda\neq 0$, we have $g_\lambda = \lambda^m g.$\\

\indent Let ${\mathbb{C}}_-$ denote the complex halfplane $\{{\rm Im}\tau < 0\}$ with closure $\overline{{\mathbb{C}}_-}$. Then:\\

\noindent {\bf {Lemma 2.3}} (\cite{BeGS}, Prop. 1.9). {\it Let $q(\xi, \tau)\in C^{\infty}(({\mathbb{R}}^n\times {\mathbb{R}})/0)$ be a parabolic homogeneous symbol
of degree $m$ such that:}\\
(i){\it ~ $q$ extends to a continuous function on $(
{\mathbb{R}}^n\times \overline{{\mathbb{C}}_-})/0 $ in such way to
be holomorphic in the
last variable when the latter is restricted to ${{\mathbb{C}}}_-$.}\\
{\it Then there is a unique $g\in  {\textsl{S}} ({\mathbb{R}}^{n+1})$ agreeing with q on ${\mathbb{R}}^{n+1}/0$ so that:}\\
(ii) {\it $g$ is homogeneous of degree $m$};\\
(iii) {\it The inverse Fourier transform $\breve{g}(x, t)$ vanishes for $t < 0.$}\\

\indent Let $U$ be an open subset of ${\mathbb{R}}^n$. We define Volterra symbols and Volterra $\Psi DO$¡¯s on $U\times {\mathbb{R}}^{n+1}/0$
as follows.\\

\noindent {\bf {Definition 2.4}}~~{\it $S_V^m(U\times
{\mathbb{R}}^{n+1}),~m\in{\mathbb{Z}}$ ,
consists of smooth functions $q(x, \xi, \tau)$ on $U\times
{\mathbb{R}}^{n}\times {\mathbb{R}}$ with
an asymptotic expansion $q\sim \sum_{j\geq 0}q_{m-j},$ where:}\\
-{\it  $q_l\in C^{\infty}(U\times [({\mathbb{R}}^{n}\times
{\mathbb{R}})/0]$ is a homogeneous Volterra symbol of degree $l$,
i.e. $q_l$ is parabolic homogeneous of degree $l$ and satisfies the
property (i) in Lemma 2.3 with respect to the last $n + 1$
variables;}\\
- {\it The sign $\sim$ means that, for any integer $N$ and any compact $K,~ U,$ there is a constant
$C_{NK\alpha\beta k}>0$ such that for $x\in K$ and for $|\xi|+|\tau|^{\frac{1}{2}}>1$ we have}

$$|\partial_x^\alpha\partial_\xi^\beta\partial_\tau^k(q-\sum_{j<N}q_{m-j})(x,\xi,\tau)|\leq
C_{NK\alpha\beta k}(|\xi|+|\tau|^{\frac{1}{2}})^{m-N-|\beta|-2k}.\eqno(2.9)$$\\

\noindent {\bf Definition 2.5}  {\it $\Psi_V^m(U\times
{\mathbb{R}}), ~m\in{\mathbb{Z}}$ , consists of
continuous operators $Q$ from $C_c^{\infty}(U_x\times
{\mathbb{R}}_t)$ to $C^{\infty}(U_x\times
{\mathbb{R}}_t)$
such that:}\\
(i) {\it $ Q $ has the Volterra property;}\\
(ii) {\it $Q = q(x,D_x,D_t) + R$ for some symbol $q$ in
$S_V^m(U\times {\mathbb{R}})$
and some smooth operator $R$.}\\

\indent In the sequel if $Q$ is a Volterra $\Psi DO$, we let $K_Q(x, y, t-s)$ denote its distribution kernel, so that
the distribution $K_Q(x, y, t)$ vanishes for $t< 0$.\\

\noindent {{\bf Definition 2.6} {\it  Let $q_m(x, \xi, \tau)\in
C^{\infty}(U\times ({\mathbb{R}}^{n+1}/0))$ be a
homogeneous Volterra symbol of order $m$ and let $g_m \in
C^{\infty}(U)\otimes {\mathbb{S}}'({\mathbb{R}}^{n+1})$ denote its unique homogeneous extension given by Lemma 2.3.
Then:}\\
- {\it $\breve{q}_m(x, y, t)$ is the inverse Fourier transform of $g_m(x, \xi, \tau )$ in the last $n + 1$ variables;}\\
- {\it $q_m(x,D_x,D_t)$ is the operator with kernel $\breve{q}_m(x, y-x, t).$}\\

\noindent{\bf Proposition 2.7}~~{\it The following properties hold.\\
1) Composition. Let $Q_j\in \Psi_V^{m_j}(U\times {\mathbb{R}}), ~j=1,2$
have symbols $q_j$ and suppose that $Q_1$ or $Q_2$ is
properly supported. Then $Q_1Q_2$ is a Volterra $\Psi DO$ of order $m_1+m_2$ with symbol $q_1\circ q_2\sim \sum\frac{1}{\alpha!}\partial^\alpha_\xi
q_1D^\alpha_xq_2.$}\\
2) {\it Parametrices. An operator $Q$ is the order $m$ Volterra
$\Psi DO$ with the paramatrix $P$ then
$$QP=1-R_1,~~~PQ=1-R_2\eqno(2.10)$$
where $R_1,~R_2$ are smooth operators.}\\

\noindent{\bf Proposition 2.8}~{\it The differential operator $D^2 +
\partial_t$ is invertible and
its inverse $(D^2 + \partial_t)^{-1}$ is a Volterra $\Psi DO$ of order $-2$.}\\

\indent   We denote by $M^\phi$ the fixed-point set of $\phi$, and for $a = 0,\cdots ,n,$ we let
   $M^\phi=\bigcup _{0\leq a\leq n} M_{a}^\phi$, where  $M_{a}^\phi$ is an $a$-dimensional submanifold. Given a fixed-point $x_0$ in a component
   $M_{a}^\phi$, consider some local coordinates $x = (x^1,\cdots , x^a)$ around
$x_0.$ Setting $b = n-a,$ we may further assume that over the range
of the domain of the local coordinates there is an orthonormal frame
$e_1(x),\cdots , e_b(x)$ of $N^\phi_z$. This defines fiber
coordinates $v = (v_1, \cdots , v_b).$ Composing with the map
$(x,v)\in N^\phi(\varepsilon_0)\rightarrow {\rm exp}_x(v)$ we then
get local coordinates $x^1,\cdots,x^a,v^1,\cdots,v^b$ for $M_z$ near
the fixed point $x_0$. We shall refer to this type of coordinates as
{\it tubular coordinates.} Then $N^\phi(\varepsilon_0)$ is homeomorphic with a tubular neighborhood of $M^\phi$. \\
\indent By the Mckean-Singer formula, we have,
$${\rm Ind}_\phi(D^+)={\rm Str}[\widetilde{\phi}e^{-tD^2}]
=\int_{M}{\rm Str}[ \widetilde{\phi}k_t(x,\phi(x))]dx
=\int_{M}{\rm Str}[
\widetilde{\phi}K_{(D^2+\partial_t)^{-1}}(x,\phi(x),t)]dx.
\eqno(2.11)$$
 Let $Q=(D^2+\partial_t)^{-1}$. For $x\in M^\phi$ and
$t>0$ set
$$I_Q(x,t):=\widetilde{\phi}(x)^{-1}\int_{N_x^\phi(\varepsilon)}\phi({\rm
exp}_xv)K_Q({\rm exp}_xv,{\rm exp}_x(\phi'(x)v),t)dv.\eqno(2.12)$$
Here we use the trivialization of $\wedge(T^*M)$ about the tubular
coordinates. Using the tubular coordinates, then
$$I_Q(x,t)=\int_{|v|<\varepsilon}\widetilde{\phi}(x,0)^{-1}\widetilde{\phi}(x,v)K_Q(x,v;x,\phi'(x)v;t)dv.\eqno(2.13)$$
Let
$$q^{\wedge(T^*M)}_{m-j}(x,v;\xi,\nu;\tau):=\widetilde{\phi}(x,0)^{-1}\widetilde{\phi}(x,v)q_{m-j}(x,v;\xi,\nu;\tau).\eqno(2.14)$$
Recall\\

\noindent {\bf Proposition 2.9} (\cite{PW}, Proposition 3.4) {\it Let $Q\in
\Psi_V^m(M\times {\mathbb{R}},\wedge(T^*M)),~m\in {\mathbb{Z}}.$ Uniformly on
each component $M_{a}^\phi$
$$I_Q(x,t)\sim \sum_{j \geq
0}t^{-(\frac{a}{2}+[\frac{m}{2}]+1)+j}I_Q^j(x) ~~~~~~{\rm
as}~~t\rightarrow 0^+,\eqno(2.15)$$ where $I_Q^j(x)$ is defined by}
$$I_Q^{(j)}(x):=\sum_{|\alpha|\leq
m-[\frac{m}{2}]+2j}\int\frac{v^\alpha}{\alpha!}\left(\partial_v^\alpha
q^{\wedge(T^*M)}_{2[\frac{m}{2}]-2j+|\alpha|}\right)^\vee(x,0;0,(1-\phi'(x))v;1)dv.\eqno(2.16)$$\\

Let
$$e(TM^\phi,\nabla^{TM^\phi})=Pf\left[-\frac{1}{2\pi}R^{TM^\phi}\right]$$
be the Euler form of $TM^\phi$ associated with $\nabla^{TM^\phi}$,
where $\nabla^{TM^\phi}$ is the Levi-Civita connection on $M^\phi$ and $R^{TM^\phi}$ its curvature. Then we have\\

\noindent {\bf Theorem 2.10} (The equivariant Gauss-Bonnet-Chern theorem) {\it The following formula holds}\\

$${\rm Ind}_\phi(D^+)=\int_{M^\phi}e(TM^\phi,\nabla^{TM^\phi}).\eqno(2.17)$$\\

\indent Let $(V,q)$ be a finite dimensional real vector space
equipped with a quadratic form. Let $C(V, q)$ be the associated
Clifford algebra, i.e. the associative algebra generated by V with
the relations $$v \cdot w + w \cdot v = -2q(v,w)$$ for $v,w \in V .$
Let $e_1\cdots, e_n$ be the orthomormal basis of $(V,q)$, Let
$C(V,q)\widehat{\otimes}C(V,-q)$ be the grading tensor product of
$C(V,q)$ and $C(V,-q)$ and $\wedge^*V\widehat{\otimes}\wedge^*V$ be
the grading tensor product of $\wedge^*V $ and $\wedge^*V$. Define
the symbol map:
$$\sigma: C(V,q)\widehat{\otimes}C(V,-q)\rightarrow \wedge^*V\widehat{\otimes}\wedge^*V;$$ $$~~\sigma(c(e_{j_1})\cdots c(e_{j_l})\otimes 1)=e^{j_1}\wedge\cdots
\wedge e^{j_l}\otimes 1;~~\sigma(1\otimes \widehat{c}(e_{j_1})\cdots
\widehat{c}(e_{j_l}))=1\otimes \widehat{e}^{j_1}\wedge\cdots \wedge
\widehat{e}^{j_l}.\eqno(2.18)$$ Using the interior multiplication
$\iota(e_j):\wedge^*V\rightarrow \wedge^{*-1}V$ and the exterior
multiplication $\varepsilon(e_j):\wedge^*V\rightarrow
\wedge^{*+1}V$, we define representations of $C(V,q )$ and $C(V,-q
)$ on the exterior algebra:
$$c: C(V,q)\rightarrow {\rm End}\wedge V,~e_j\mapsto c(e_j):\varepsilon(e_j)-\iota(e_j),$$
$$\widehat{c}: C(V,-q)\rightarrow {\rm End}\wedge V,~e_j\mapsto \widehat{c}(e_j):\varepsilon(e_j)+\iota(e_j),$$
The tensor product of these representations yields an isomorphism of
superalgebras
$$c\otimes \widehat{c}: C(V,q)\widehat{\otimes}C(V,-q)\rightarrow {\rm End}\wedge V,$$
which we will also denote by $c$. We obtain a supertrace (i.e. a linear functional
vanishing on supercommutators) on $C(V,q)\widehat{\otimes}C(V,-q)$ by setting
${\rm Str}(a) = {\rm Str}_{{\rm End}\wedge V }[c(a)]$
for $a \in C(V,q)\widehat{\otimes}C(V,-q)$, where ${\rm Str}_{{\rm End}\wedge V }$ is the canonical supertrace on ${\rm End}V$.\\

\noindent{\bf Lemma 2.11} {\it For $1\leq i_1<\cdots <i_p\leq n,~1\leq j_1<\cdots <j_q\leq n,$ we have}
$${\rm Str}[c(e_{i_1})\cdots c(e_{i_p})\widehat{c}({e}_{j_1})\cdots \widehat{c}({e}_{j_q})]
=(-1)^{\frac{n}{2}}2^n,$$ {\it when $p=q=n$, otherwise equals zero.}\\

We will also denote the volume element in $\wedge V\widehat{\otimes}\wedge V$ by $\omega=e^1\wedge\cdots\wedge e^n\wedge \widehat{e}^1\wedge\cdots\wedge\widehat{e}^n$ and For $a\in \wedge V\widehat{\otimes}\wedge V$ let $Ta$ be the coefficient of $\omega$ in $a$.
The linear functional
$T : \wedge V\widehat{\otimes}\wedge V\rightarrow R$ is called the Berezin trace. Then
for $a \in C(V,q)\widehat{\otimes}C(V,-q)$ , one has ${\rm Str}(a) = (-1)^{\frac{n}{2}}2^n(T\sigma)(a).$
 We define the Getzler order as follows:
 $${\rm deg}\partial_j=\frac{1}{2}{\rm deg}\partial_t=2{\rm deg}c(e_j)=2{\rm deg}\widehat{c}(e_j)=-{\rm deg}x^j=1.\eqno(2.19)$$
 Let $Q\in \Psi_V^*({\mathbb{R}}^n\times {\mathbb{R}}, \wedge ^*T^*M)$ have symbol $$q(x,\xi,\tau)\sim
 \sum_{k\leq m'}q_{k}(x,\xi,\tau),\eqno(2.20)$$
 where $q_{k}(x,\xi,\tau)$ is an order $k$ symbol. Then taking components
in each subspace $\wedge^jT^*M$
and using Taylor expansions at $x = 0$ give formal expansions
$$\sigma[q(x,\xi,\tau)]\sim\sum_{j,k}\sigma[q_{k}(x,\xi,\tau)]^{(j)}\sim\sum_{j,k,\alpha}\frac{x^\alpha}{\alpha!}
\sigma[\partial_x^\alpha
q_{k}(0,\xi,\tau)]^{(j)}.\eqno(2.21)$$ The symbol
$\frac{x^\alpha}{\alpha!} \sigma[\partial_x^\alpha
q_{k}(0,\xi,\tau)]^{(j)}$ is the Getzler homogeneous
of degree $k+\frac{j}{2}-|\alpha|$. So we can expand $\sigma[q(x,\xi,\tau)]$ as
$$\sigma[q(x,\xi,\tau)]\sim \sum_{j\geq 0}q_{(m-\frac{j}{2})}(x,\xi,\tau),~~~~~~~~~q_{(m)}\neq 0, \eqno(2.22)$$
where $q_{(m-\frac{j}{2})}$ is a Getzler homogeneous symbol of degree $m-\frac{j}{2}$.\\

\noindent {\bf Definition 2.12} The $m$ is called as the Getzler order of $Q$. The symbol $q_{(m)}$ is the principle Getzler
homogeneous symbol of $Q$. The operator $Q_{(m)}=q_{(m)}(x,D_x,D_t)$ is called as the model operator of $Q$.\\

\indent Let $e_1, \dots , e_n$ be an oriented orthonormal basis of
$T_{x_0}M$ such that $e_1,\cdots , e_a$ span $T_{x_0}M^\phi$ and
$e_{a+1},\cdots , e_n$ span $N_{ x_0}^\phi$ . This provides us with
normal coordinates $(x_1, \cdots , x_n)\rightarrow {\rm
exp}_{x_0}(x^1e_1+\cdots+x^ne_n).$ Moreover using parallel
translation enables us to construct a synchronous local oriented
tangent frame $e_1(x), . . . , e_n(x)$ such that $e_1(x),\cdots ,
e_a(x)$ form an oriented frame of $TM_{z,a}^\phi$ and $e_{a+1}(x),
\cdots , e_n(x)$ form an (oriented) frame $N^\phi$ (when both frames
are restricted to $M^\phi).$ This gives rise to trivializations of
the tangent and exterior algebra bundles. Write
$$\phi'(0)=\left(\begin{array}{lcr}
  1  & 0  \\
   0  &  \phi^N
\end{array}\right).
$$
Let $\wedge(n)=\wedge^*{\mathbb{R}}^n$ be the
exterior algebra of ${\mathbb{R}}^n$. We shall use the
following gradings on $\wedge(n)\widehat{\otimes}\wedge (n),$
$$\wedge(n)\widehat{\otimes}\wedge (n)=\bigoplus_{\begin{array}{lcr}
  1\leq k_1,k_2\leq a \\
  1\leq \overline{l_1},\overline{l_2} \leq b
\end{array}}\wedge^{k_1,\overline{l_1}}(n)\widehat{\otimes}\wedge^{k_2,\overline{l_2}}(n),$$
 where
 $\wedge^{k,\overline{l}}(n)$
is the space of forms $dx^{i_1}\wedge\cdots\wedge
dx^{i_{k+\overline{l}}}$ with $1\leq i_1<\cdots <i_k\leq a$ and $a +
1\leq i_{k+1} < \cdots < i_{k+\overline{l}}\leq n.$ Given a form
$\omega\in\wedge (n)\widehat{\otimes}\wedge (n)$ we shall denote by
$\omega^{((k_1,\overline{l_1}),(k_2,\overline{l_2}))}$ its component in $\wedge^{k_1,\overline{l_1}}(n)\widehat{\otimes}\wedge^{k_2,\overline{l_2}}(n).$
 We denote by
$|\omega|^{((a,0),(a,0))}$ the Berezin integral $|\omega^{((*,0),(*,0))}|^{((a,0),(a,0))}$
of its component $\omega^{((*,0),(*,0))}$ in $\wedge^{((*,0),(*,0))}(n).$
Similar to Lemma 3.6 in \cite{Wa}. we have by (2.19)\\

\noindent{\bf Lemma 2.13}~ {\it $Q\in \Psi_V^*({\mathbb{R}}^n\times
{\mathbb{R}},\wedge(T^*M))$ has the Getzler order $m$
and model operator $Q_{(m)}$. Let $j$ be even, then as $t\rightarrow 0^+$}\\
\indent (1) $\sigma[I_Q(0,t)]^{(j)}=O(t^{\frac{\frac{j}{2}-m-a-1}{2}})$ {\it if
$m-\frac{j}{2}$
is odd.}\\
\indent (2)
$\sigma[I_Q(0,t)]^{(j)}=O(t^{\frac{\frac{j}{2}-m-a-2}{2}})I_{Q(m)}(0,1)^{(j)}+O(t^{\frac{\frac{j}{2}-m-a}{2}})$ {\it if
$m-\frac{j}{2}$}
is even.}\\
{\it In particular, for $m=-2$, $j=2a$ and $a$ is even we get}
$$\sigma[\psi_tI_Q(0,t)]^{((a,0),(a,0))}=I_{Q(-2)}(0,1)^{((a,0),(a,0))}+O(t^{\frac{1}{2}}).\eqno(2.23)$$\\

\indent By the Weitzenb\"{o}ck formula, we have
$$D^2=-\sum_{j=1}^n(\nabla_{e_j}^2-\nabla_{\nabla^{TM}_{e_j}e_j})+\frac{r_M}{4}-\frac{1}{8}\sum_{1\leq i,j,k,l\leq n}R_{ijkl}c(e_i)c(e_j)\widehat{c}(e_k)\widehat{c}(e_l).\eqno(2.24)$$
By (2.19) and (2.24), we get the model operator of $\frac{\partial}{\partial t}+D^2$ is
$$\frac{\partial}{\partial t}-\sum_{j=1}^n\frac{\partial^2}{\partial y^2_j}
-\frac{1}{8}\sum_{1\leq i,j,k,l\leq n}R_{ijkl}e^i\wedge e^j\wedge\widehat{e^k}\wedge\widehat{e^l}.$$
By
$$(\frac{\partial}{\partial t}-\sum_{j=1}^n\frac{\partial^2}{\partial y^2_j}
-\frac{1}{2}\dot{R})K_{Q(-2)}(x,y,t)=0,\eqno(2.25)$$
we get
$$K_{Q(-2)}(x,y,t)=(4\pi t)^{\frac{n}{2}}{\rm exp}(-\frac{1}{4t}||x-y||^2)e^{\frac{t\dot{R}}{2}}.\eqno(2.26)$$
Similar to Lemma 9.13 in \cite{PW}, we get
$$I_{Q(-2)}(0,t)=(4\pi t)^{-\frac{a}{2}}{\rm det}^{-1}(1-\phi^N)e^{\frac{t\dot{R}}{2}}.\eqno(2.27)$$
Let the matrix $\phi^N$ equal
$$\phi^N={\rm diag}\left(\left(\begin{array}{lcr}
 {\rm cos}\theta_{{\frac{a}{2}+1 }} & {\rm sin}\theta_{{\frac{a}{2}+1 }}  \\
   -{\rm sin}\theta_{{\frac{a}{2}+1 }}  &  {\rm cos}\theta_{{\frac{a}{2}+1 } }
\end{array}\right),\cdots ,\left(\begin{array}{lcr}
  {\rm cos}\theta_{\frac{n}{2} } & {\rm sin}\theta_{\frac{n}{2} }  \\
   -{\rm sin}\theta_{\frac{n}{2} }  &  {\rm cos}\theta_{\frac{n}{2} }
\end{array}\right)\right).$$
We have the lemma\\

\noindent{\bf Lemma 2.14} (\cite{Zh}, Lemma 3.2)~{\it The following equality holds}\\
$$\widetilde{\phi}=(\frac{1}{2})^{\frac{n-a}{2}}\prod_{j=\frac{a}{2}+1}^{\frac{n}{2}}\left[(1+{\rm cos}\theta_j)-(1-{\rm cos}\theta_j)c(e_{2j-1})
c(e_{2j})\hat{c}(e_{2j-1})
\widehat{c}(e_{2j})\right.$$
$$\left.+{\rm sin}\theta_j(c(e_{2j-1})
c(e_{2j})-\hat{c}(e_{2j-1})
\widehat{c}(e_{2j}))\right].\eqno(2.28)$$\\

By Lemma 2.14, we have
$$\sigma[\widetilde{\phi}]^{((0,b),(0,b))}=(\frac{1}{2})^{\frac{n-a}{2}}(-1)^{\frac{b}{2}}\prod_{j=\frac{a}{2}+1}^{\frac{n}{2}}
(1-{\rm cos}\theta_j)e^{a+1}\wedge\cdots\wedge e^n\wedge \widehat{e}^{a+1}\wedge\cdots\wedge \widehat{e}^n$$
$$=(-\frac{1}{4})^{\frac{b}{2}}{\rm det}(1-\phi^N)e^{a+1}\wedge\cdots\wedge e^n\wedge \widehat{e}^{a+1}\wedge\cdots\wedge \widehat{e}^n.\eqno(2.29)$$
So we get the following lemma\\

\noindent{\bf Lemma 2.15} ~{\it For $A\in C(V,q)\widehat{\otimes}C(V,-q)$, we have}
$${\rm Str}[\widetilde{\phi}A]=(-1)^{\frac{n}{2}}2^n(-\frac{1}{4})^{\frac{b}{2}}{\rm det}(1-\phi^N)
|\sigma(A)|^{((a,0),(a,0))}$$
$$+(-1)^{\frac{n}{2}}2^n\sum_{l_1,{\rm {or}}~ l_2<b}|\sigma(\widetilde{\phi})^{((0,l_1),(0,l_2))}
\sigma(A)^{((a,b-l_1),(a,b-l_2))}|^{(n,n)}.\eqno(2.30)$$\\

By (2.27), (2.30) and Lemma 2.13, we get
$${\rm lim}_{t\rightarrow 0}{\rm
Str}[\widetilde{\phi}I_{(D^2+\partial_t)^{-1}}(x_0,t)]
=(-1)^{\frac{n}{2}}2^n(-\frac{1}{4})^{\frac{b}{2}}(4\pi
)^{-\frac{a}{2}}|e^{\frac{\dot{R}}{2}}|^{((a,0),(a,0))}=e(TM^\phi).\eqno(2.31)$$
So we get Theorem 2.10.\\

\section{The variation formulas for the equivariant Ray-Singer
metric}
   \quad Let $M$ be an even dimensional $n$ oriented closed manifold and $(F,\nabla^F)$ be a flat complex vector bundle over $M$. Let
    $g^{TM}$ be a Riemannian metric on $M$ and $h^F$ be a Hermitian metric on $F$. We will not assume $h^F$ to be parallel with respect to
    $\nabla^F.$ Let $G$ be a compact Lie group acting smoothly on $M$ such that the metric $g^{TM},~h^F$ and the flat connection $\nabla^F$ are
    preserved. Let $A^*(M,F)=\Gamma(M,\wedge T^*M\otimes F)$ denote
    the differential forms on $M$ with values in $F$. Let
    $d^F:A^*(M,F)\rightarrow A^{*+1}(M,F)$ denote the exterior
    differential associated with the flat connection $\nabla^F$. The
    Hodge Laplacian is given by
    $$\triangle^F=d^Fd^{F,*}+d^{F,*}d^{F},$$ where $d^{F,*}$ denotes
    the formal adjoint of $d^{F}$. For $t>0$ let ${\rm
    exp}(-t\triangle^F)$ denote the heat operator. We consider
    $1$-parameter families of $G$-invariant metrics: \\
    \indent (1)
    $\varepsilon\mapsto g^{TM}(\varepsilon)$ with $g^{TM}(0)=g^{TM}$~ (2) $\varepsilon\mapsto h^{F}(\varepsilon)$ with
    $h^{F}(0)=h^{F}.$ \\ Let $C=\star^{-1}\dot{\star}$ and
    $V=(h^F)^{-1}\dot{h^F}$. For a $1$-parameter families of Riemannian
    metrics $\varepsilon\mapsto g^{TM}(\varepsilon)$ we set
    $$\dot{S}:=\dot{\nabla}^{TM}-\frac{1}{2}[\nabla^{TM},(g^{TM})^{-1}\dot{g}^{TM}]\in
    A(M,so(TM)).$$
    We define the transgression form
    $$\widetilde{e}'(TM):=\frac{\partial}{\partial b}|_{b=0}{\rm
    Pf}\left[-\frac{1}{2\pi}(R^{TM}+b\dot{S})\right].$$
Since $h^F$ is not necessarily parallel with respect to $\nabla^F$, we may
define a second flat connection $(\nabla^F)^T$ on $F$ by the formula
$$(\nabla^F)^T=(h^F)^{-1}\nabla^{F,*}h^F,$$ where $\nabla^{F,*}$
denotes the connection induced by $\nabla^{F}$ on $F^*$ and
$h^F:F\rightarrow F^*$ is the isomorphism induced by $h^F$. Observe
that $(\nabla^F)^T=\nabla^F$ if and only if $\nabla^Fh^F=0$. We set
$$\omega(F,h^F):=(\nabla^F)^T-\nabla^F\in A^1(M.{\rm End}(F))$$ and for $\phi\in G$~
$$\theta(\phi,F,h^F):={\rm Tr}[\phi
\omega(F,h^F)]\in A^1(M^\phi).$$ Then $\theta(\phi,F,h^F)$ is closed
and that its cohomology class does not depend on $h^F$. In the
following, we will prove\\

\noindent {\bf Theorem 3.1 (Bismut-Zhang)} {\it For $\phi\in G$, we
have}
$${\rm lim}_{t \rightarrow 0}{\rm Str}[\widetilde{\phi}V{\rm
exp}(-t\triangle^F)]=\int_{M^\phi}{\rm
Tr}[\phi^FV]e(TM^\phi,\nabla^{TM^\phi}),\eqno(3.1)$$
$${\rm lim}_{t \rightarrow 0}{\rm Str}[\widetilde{\phi}C{\rm
exp}(-t\triangle^F)]=-\int_{M^\phi}\theta(\phi,F,h^F)\widetilde{e}'(TM^\phi),\eqno(3.2)$$

\indent By Theorem 3.1, we can easily get the variation formulas for the
equivariant Ray-Singer metric due to Bismut-Zhang (see pp.4-5 in \cite{We}).
In general, neither of the connections $\nabla^F$ and $(\nabla^F)^T$ will preserve the metric
$h^F.$ As in \cite{BZ1} we define a third connection on $F$
 $$\nabla^{F,e}=\frac{1}{2}(\nabla^F+(\nabla^F)^T).$$
 This connection will preserve $h^F$, but it will in general not be flat.\\
\indent In the following we will write ${\mathcal{E}}=\wedge
T^*M\otimes F$. We will also denote by $\nabla^{F,e}$ the tensor
product connection $$\nabla^{\wedge T^*M}\otimes 1+1\otimes
\nabla^{F,e}$$ on ${\mathcal{E}},$ where $\nabla^{\wedge T^*M}$ is
the connection on ${\wedge T^*M}$ induced by $\nabla^{TM}$. Let
$\triangle^{{\mathcal{E}},e}$ denote the connection Laplacian on
${\mathcal{E}}$ associated to the connection $\nabla^{F,e}$. Since
$\nabla^{F,e}$ is a metric connection on ${\mathcal{E}}$, the
operator $\triangle^{{\mathcal{E}},e}$ will be formally
self-adjoint.\\

\noindent{\bf Proposition 3.2 }(\cite{BZ1}, Theorem 4.13): (Lichnerowicz formula for $\triangle(F)$)~ {\it One has
$$ \triangle(F)=-\triangle^{{\mathcal{E}},e}+E$$
with $E\in \Gamma(M,{\rm End} {\mathcal{E}})$ which w.r.t. a local ON-frame ${e_j}$ is given by}
\begin{eqnarray*}
E&=&-\frac{1}{8}\sum_{i,j,k,l}(R^{TM}(e_i,e_j)e_k,e_l)c(e_i)c(e_j)\widehat{c}(e_k)\widehat{c}(e_l)\\
& ~ &-\frac{1}{8}\sum_{i,j}c(e_i)c(e_j)\omega(F,h^F)^2(e_i,e_j)+\frac{1}{8}\sum_{i,j}\widehat{c}(e_i)\widehat{c}(e_j)\omega(F,h^F)^2(e_i,e_j)\\
& ~ & -\frac{1}{2}\sum_{i,j}c(e_i)\widehat{c}(e_j)[\nabla^{T^*M\otimes {\rm End }F}\omega(F,h^F)(e_j)+
\frac{1}{2}\omega(F,h^F)^2(e_i,e_j)]\\
& ~ &
+\frac{1}{4}\sum_j(\omega(F,h^F)(e_j))^2+\frac{1}{4}r^M,
\end{eqnarray*}
{\it where $r^M$ denotes the scalar curvature of $(M, g^{TM})$.}\\

Let $$D^0=\sum_{j=1}^nc(e_j)\nabla^e_{e_j}$$ be a self-adjoint twisted Dirac operator. Then
$$(D^0)^2=-\triangle^{{\mathcal{E}},e}-\frac{1}{8}\sum_{i,j,k,l}(R^{TM}(e_i,e_j)e_k,e_l)c(e_i)c(e_j)\widehat{c}(e_k)\widehat{c}(e_l)$$
$$-\frac{1}{8}\sum_{i,j}c(e_i)c(e_j)\omega(F,h^F)^2(e_i,e_j)+\frac{1}{4}r^M.\eqno(3.3)$$
Let
$$L(\omega)=\frac{1}{8}\sum_{i,j}\widehat{c}(e_i)\widehat{c}(e_j)\omega(F,h^F)^2(e_i,e_j)
-\frac{1}{2}\sum_{i,j}c(e_i)\widehat{c}(e_j)[\nabla^{T^*M\otimes {\rm End }F}\omega(F,h^F)(e_j)$$
$$+
\frac{1}{2}\omega(F,h^F)^2(e_i,e_j)]
+\frac{1}{4}\sum_j(\omega(F,h^F)(e_j))^2.\eqno(3.4)$$
Then $$\triangle(F)=(D^0)^2+L(\omega).\eqno(3.5)$$

\noindent {\bf The Proof of (3.1):}\\

\indent By (2.19), we get the model operator of $\frac{\partial}{\partial t}+\triangle(F)$ is still
$$\frac{\partial}{\partial t}-\sum_{j=1}^n\frac{\partial^2}{\partial y^2_j}
-\frac{1}{2}\dot{R}.$$ Let $P$ be a differential operator, then we have
$${\rm Str}[\widetilde{\phi}P{\rm exp}(-t\triangle(F))]=\int_{M^\phi}{\rm Str}[\widetilde{\phi}I_{P(\frac{\partial}{\partial t}+\triangle(F))^{-1}}(x,t)]dx+O(t^{\infty}).\eqno(3.6)$$
When $P=V$ is a $0$-order differential operator and $O_G(V)=0$, we have $$O_G(V(\frac{\partial}{\partial t}+\triangle(F))^{-1})=-2,$$
and the model operator of $V(\frac{\partial}{\partial t}+\triangle(F))$ is
$$V\left[\frac{\partial}{\partial t}-\sum_{j=1}^n\frac{\partial^2}{\partial y^2_j}
-\frac{1}{2}\dot{R}\right].$$
Then similar to the proof of Theorem 2.10, we get (3.1). ~~~$\Box$\\

\noindent {\bf Lemma 3.3} (\cite{BZ1}, Proposition 4.15) {\it The endomorphism $C=\star^{-1}\dot{\star}$ is given in terms of Clifford variable by}
$$C=-\frac{1}{2}\sum_{i,j}((g^{TM})^{-1}\dot{g}^{TM}e_i,e_j)c(e_i)\widehat{c}(e_j),\eqno(3.7)$$
$$\sigma(C)=-\frac{1}{2}\sum_{i,j}((g^{TM})^{-1}\dot{g}^{TM}e_i,e_j)e^i\wedge\widehat{e}^j.\eqno(3.8)$$\\

As in \cite{BZ1}, we introduce an auxiliary even Clifford variable, such that $\sigma^2=1$ and $\sigma$ commutes with all the previously considered
operators. Let $A,B\in{\rm End}({\mathcal{E}})$ be trace classes. Then $A+\sigma B$ lies in ${\rm End}({\mathcal{E}})\otimes {\mathbb{R}}(\sigma)$. Set
$${\rm Str}^\sigma[A+\sigma B]={\rm Str}[B].\eqno(3.9)$$
Let
$$(D^2)^{\rm odd}:=-\frac{1}{2}\sum_{i,j}c(e_i)\widehat{c}(e_j)[\nabla^{T^*M\otimes {\rm End }F}\omega(F,h^F)(e_j)+
\frac{1}{2}\omega(F,h^F)^2(e_i,e_j)],\eqno(3.10)$$
$$(D^2)^{\rm even}:=\triangle(F)-(D^2)^{\rm odd}.\eqno(3.11)$$
By Theorem 4.18 in \cite{BZ2}, we have
$${\rm Str}[\widetilde{\phi}C{\rm exp}(-t\triangle(F))]={\rm Str}^\sigma[\widetilde{\phi}C{\rm exp}(-t((D^2)^{\rm even}
+\sigma (D^2)^{\rm odd}))].\eqno(3.12)$$
Let
$$L(\omega,\sigma)=\frac{1}{8}\sum_{i,j}\widehat{c}(e_i)\widehat{c}(e_j)\omega(F,h^F)^2(e_i,e_j)
-\frac{1}{2}\sigma\sum_{i,j}c(e_i)\widehat{c}(e_j)[\nabla^{T^*M\otimes {\rm End }F}\omega(F,h^F)(e_j)$$
$$+
\frac{1}{2}\omega(F,h^F)^2(e_i,e_j)]
+\frac{1}{4}\sum_j(\omega(F,h^F)(e_j))^2.\eqno(3.13)$$
Then $$(D^2)^{\rm even}
+\sigma (D^2)^{\rm odd}=(D^0)^2+L(\omega,\sigma).\eqno(3.14)$$
By the Duhamel principle, we have, when $t$ is small,
$${\rm Str}^\sigma[\widetilde{\phi}C{\rm exp}(-t((D^0)^2+L(\omega,\sigma))]
=(\sum_{k=1}^{\frac{n}{2}}+\sum_{k>\frac{n}{2}})(-t)^k\int_{\triangle_k}{\rm
Str}^\sigma[\widetilde{\phi}C$$
$$\cdot
e^{-t_0t(D^0)^2}L(\omega,\sigma)
\cdots L(\omega,\sigma)e^{-t_kt(D^0)^2}]{\rm dvol}_{\triangle_k}.\eqno(3.15)$$
By the H$\ddot{o}$lder inequality and the Weyl theorem, we can get
$$\sum_{k>\frac{n}{2}}t^k\int_{\triangle_k}\left|{\rm Str}^\sigma[\widetilde{\phi}C
\cdot
e^{-t_0t(D^0)^2}L(\omega,\sigma)
\cdots L(\omega,\sigma)e^{-t_kt(D^0)^2}]\right|{\rm dvol}_{\triangle_k}$$
$$\leq \sum_{k>\frac{n}{2}}\frac{t^k||L(\omega,\sigma)||^k}{k!}O(t^{-\frac{n}{2}})=O(t^{\frac{1}{2}}).\eqno(3.16)$$
\indent Let ${B}$ be an
operator and $l$ be a positive interger.
 Write
$${B}^{[l]}=[(D^0)^2,{B}^{[l-1]}],~{B}^{[0]}={B}.$$\\

\noindent{\bf Lemma 3.4} (\cite{Fe}, \cite{CH}, Lemma 2.7])~{\it Let ${B}$ a finite order
differential
operator. Then for any $s>0$, we have:}\\
$$e^{-s(D^0)^2}{B}=\sum^{N-1}_{l=0}\frac{(-1)^l}{l!}s^l{B}^{[l]}e^{-s(D^0)^2}+(-1)^Ns^N{B}^{[N]}(s),\eqno(3.17)$$
\noindent {\it where ${B}^{[N]}(s)$ is given by}\\
$${B}^{[N]}(s)=\int_{\triangle_N}e^{-u_1s(D^0)^2}{B}^{[N]}e^{-(1-u_1)s(D^0)^2}du_1du_2\cdots du_N.\eqno(3.18)$$\\

Similar to Theorem 1 in \cite{CH}, we have
$$t^k\int_{\triangle_k}{\rm Str}^\sigma[\widetilde{\phi}C
\cdot
e^{-t_0t(D^0)^2}L(\omega,\sigma)
\cdots L(\omega,\sigma)e^{-t_kt(D^0)^2}]{\rm dvol}_{\triangle_k}$$
$$
=\sum_{0\leq |\lambda|\leq n-k}c_{\lambda,k}t^{k+\lambda_1+\cdots \lambda_k}\int_{\triangle_k}{\rm Str}^\sigma[\widetilde{\phi}C
[L(\omega,\sigma)]^{[\lambda_1]}\cdots L(\omega,\sigma)]^{[\lambda_k]}e^{-t(D^0)^2}]+O(t^{\frac{1}{2}}),\eqno(3.19)$$
where $c_{\lambda,k}$ is a constant. Since $O_G((D^0)^2)=2$ and $O_G(L(\omega,\sigma))=1$ and $O_G(C)=1$, we have
$$O_G(CL(\omega,\sigma)^{[\lambda_1]}\cdots L(\omega,\sigma)^{[\lambda_k]})=k+1+2(\lambda_1+\cdots+\lambda_k).\eqno(3.20)$$
So when $k>1$, by Lemma 2.13, we get
$$t^k\int_{\triangle_k}{\rm Str}^\sigma[\widetilde{\phi}C
\cdot
e^{-t_0t(D^0)^2}L(\omega,\sigma)
\cdots L(\omega,\sigma)e^{-t_kt(D^0)^2}]{\rm dvol}_{\triangle_k}=O(t^{\frac{1}{2}}).\eqno(3.21)$$
So we need only compute the term
$$
-\sum_{0\leq \lambda_1\leq
n-k}c_{\lambda_1,k}t^{1+\lambda_1}\int_{\triangle_1}{\rm
Str}^\sigma[\widetilde{\phi}C
[L(\omega,\sigma)]^{[\lambda_1]}e^{-t(D^0)^2}].\eqno(3.22)$$ Since
$$(D^0)^2=-\sum_{j=1}^n\frac{\partial^2}{\partial y^2_j}
-\frac{1}{2}\dot{R}+O_G(1),\eqno(3.23)$$ we can get
$O_G([(D^0)^2,L(\omega,\sigma)])=2.$ Then
$O_G(L(\omega,\sigma)^{[\lambda_1]})=2\lambda_1$ for $\lambda_1>0$.
So by Lemma 2,13, we get
$$
-\sum_{0\leq \lambda_1\leq
n-k}c_{\lambda_1,k}t^{1+\lambda_1}\int_{\triangle_1}{\rm
Str}^\sigma[\widetilde{\phi}C
[L(\omega,\sigma)]^{[\lambda_1]}e^{-t(D^0)^2}]$$
$$
=-{\rm Str}^\sigma[t\widetilde{\phi}C
L(\omega,\sigma)e^{-t(D^0)^2}]+O(t^{\frac{1}{2}})
=-{\rm Str}[t\widetilde{\phi}C
(D^2)^{\rm odd}e^{-t(D^0)^2}]+O(t^{\frac{1}{2}})
,\eqno(3.24)$$
Similar to Theorem 2.10, we can get
$$-{\rm Str}[t\widetilde{\phi}C
(D^2)^{\rm odd}e^{-t(D^0)^2}]
=-(-\frac{1}{\pi})^{\frac{a}{2}}\left|\frac{1}{2}\left\{\sum_{i,j=1}^a((g^{TM})^{-1}\dot{g}^{TM}e_i,e_j)e^i\wedge\widehat{e}^j\right\}\right.$$
$$\left.\wedge
{\rm exp}(-\frac{\dot{R}^{TM^\phi}}{2})\wedge \frac{1}{2}\left\{\sum_{i,j=1}^a e^i\wedge\widehat{e}^j{\rm Tr}[\phi^F(\nabla_{e_i}\omega
(e_j))]\right\}\right|^{(a,0),(a,0)}.\eqno(3.25)$$
By (3.25) and using the same calculations as in the non-equivariant case (see \cite{BZ1}, p.76), we can prove (3.2).\\

 \noindent {\bf Acknowledgement.} This work was supported by NSFC. 11271062 and NCET-13-0721. The author would
like to thank the referee for careful reading and helpful comments.

 \indent{\it School of Mathematics and Statistics, Northeast Normal University, Changchun Jilin, 130024, China }\\
 \indent E-mail: {\it wangy581@nenu.edu.cn}\\

\end{document}